\documentclass[12pt]{amsart}

\usepackage{amssymb}

\newtheorem{thmA}{Theorem}

\newtheorem{theorem}{Theorem}[section]
\newtheorem{thm}[theorem]{Theorem}
\newtheorem{lemma}[theorem]{Lemma}
\newtheorem{prop}[theorem]{Proposition}

\theoremstyle{definition}
\newtheorem{definition}[theorem]{Definition}

\newtheorem{example}[theorem]{Example}
\newtheorem{remark}[theorem]{Remark}
\newtheorem{para}[theorem]{}

\def\Cal{\mathcal}
\def\Bbb{\mathbb}  
\def\A{\mathcal A}

\def\C{\mathcal C}
\def\P{\mathcal P}
\def\S{\Sigma} 
\def\G{\Gamma}
\def\g{\gamma}

\def\ARp{\langle A \mid R \rangle}
\def\s{\sigma}

\def\-{\overline}

\input epsf

\newcount\figno
\global\figno=1
\def\figure#1{\midinsert
\centerline{\epsfbox{#1}}
\centerline{Figure \the\figno}
\global\advance\figno by 1
\endinsert}

\def\scaledfigure#1#2#3{
\centerline{\epsfxsize=#2\epsfbox{#1}}
\centerline{Figure \the\figno: {#3}}
\global\advance\figno by 1}

\def\<{\langle}
\def\>{\rangle}

\def\R{\Cal R}
\def\ssm{\smallsetminus}
\def\im{\text{\rm{im} }}

\catcode`\@=11
\def\serieslogo@{\relax}
\def\@setcopyright{\relax}
\catcode`\@=12

\begin{document}
\title[Decision Problems for Combable Groups]
{The conjugacy and isomorphism problems\\
for combable groups}

\author[Bridson]{Martin R.~Bridson}
\address{Martin~R.~Bridson\\
Department of Mathematics\\
Imperial College \\
Huxley Building\\
180 Queen's Gate\\
London SW7 2BZ\\
U.K. }
\email{m.bridson@ma.ic.ac.uk}

\begin{abstract} There exist combable groups in which
the conjugacy problem is unsolvable.
The isomorphism problem is  unsolvable for
certain recursive sequences of finite presentations of 
combable groups.
\end{abstract}

\thanks{The author's research is supported
by an EPSRC Advanced Fellowship}


\date{Nov 27, 2001; columbia}


\maketitle

\centerline{\large Introduction}
\smallskip

The class of combable groups 
is one of a number of classes  that have been extensively
studied in the last fifteen years in connection with
manifestations of non-positive curvature in group theory.
The other classes closely associated to it
are word-hyperbolic
groups, automatic groups, semihyperbolic groups, bicombable
groups, and the
fundamental groups of compact non-positively curved (orbi)spaces.

Each of these classes is defined in terms of a 
convexity condition, the
common core of which is the ``fellow-traveller'' condition
that forms the definition of a combable group:
a group $\G$ with finite generating set $\A$ is
said to be {\em combable} if there is a family of
words $\{\sigma_\g \mid \g\in\G\}$ in the letters $\A^{\pm 1}$
and a 
constant $k>0$ such that for each $\g\in \G$ and $a\in\A$, the paths
in the Cayley graph $\C_\A\G$ that begin at the identity vertex
and are labelled  $\s_\g$ and $\s_{\g a}$ remain
uniformly $k$-close.

\scaledfigure{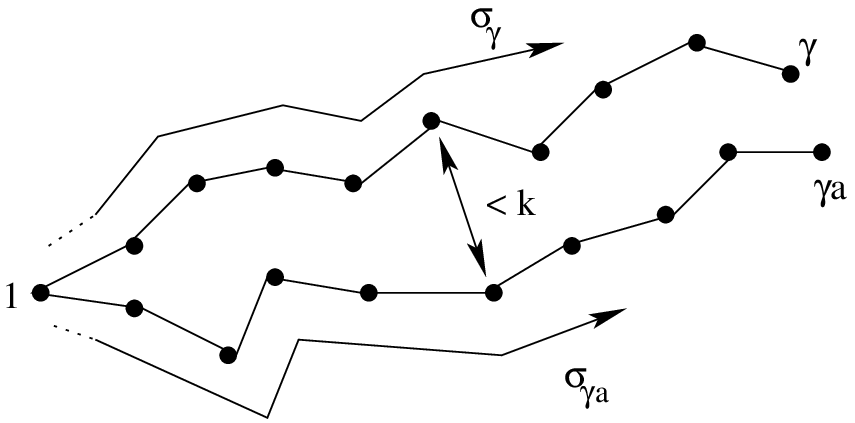}{3in}
{the fellow-traveller property}
\medskip

The results in  \cite{mb-cmh} established 
that the class of combable groups
is strictly larger than all of the other classes
listed above.  
There is an effective solution to the
 word problem in any  combable group \cite{ep+},
\cite{mb-plms}.

The conjugacy problem is known to be solvable in all of the above
classes except for the automatic groups and the combable groups, in
which cases it has been the subject of considerable speculation.

\begin{thmA} \label{thmA}
There exist combable groups in which the 
conjugacy problem is unsolvable.
\end{thmA}

The smallest and most intensively studied of the above classes 
is that formed by the hyperbolic groups. A remarkable theorem
of Zlil Sela shows that the isomorphism problem is solvable
among hyperbolic groups\footnote{In the reference cited it is
assumed that the groups are torsion-free and freely indecomposable
but Sela has a proof in the general case.}
\cite{sela}. In the light of the work of Farrell and Jones on
topological rigidity \cite{FJ}, Sela's result implies that the
homeomorphism problem is solvable among 
closed (high-dimensional) manifolds that
admit metrics of negative curvature. The results of
Farrell and Jones remain valid for
non-positively curved manifolds, and there is 
therefore considerable interest in the open
question of whether or not the isomorphism problem is solvable
among presentations of the fundamental groups of such manifolds
--- cf. \cite{BH} pp.~494. The solvability
of the isomorphism problem in the 
related classes of groups listed above has also remained open.

\begin{thmA}\label{thmB} The isomorphism problem is unsolvable in the
class of combable groups.
\end{thmA}

In order to prove this theorem we shall exhibit 
recursive sequences  of
finite subsets $\R_n$
of a finitely generated free group $F=F(\A)$  such that
each of the groups
$\G_n=\<\A \mid\R_n\>$ is combable, $n=1,2,\dots$,  but there is no algorithm
to decide which are isomorphic to $\G_1$.  
In order
to obtain such sequences we combine the 
construction of \cite{mb-cmh}
with a suitable encoding of the fact that there is no
algorithm to decide which $m$-element subsets generate
a direct product of free groups, and likewise for certain hyperbolic
groups. 

The groups $\G_n$ that we construct  are neither bicombable nor
automatic, and the isomorphism problem  remains open
for these classes.

\section{The Seeds of Undecidability}

The existence among combable groups of the
undecidability phenomenon 
asserted in Theorem \ref{thmB}  will be
established by means of a suitable encoding of the fact
that many finitely presented
groups $\G$ associated with non-positive curvature
have  the  following property:
\begin{quote}
{\bf $\Psi$: }
{\rm If $X$ is
a finite generating set for $\G$
and $m$ is a sufficiently large integer,
then there exists a recursive sequence $(S_n)$ of subsets of
the free group $F(X)$  such that each $S_n$ has
cardinality $m$ and there is no algorithm
to determine whether or not $\langle S_n\rangle = \G$; 
moreover if
equality fails then
$\langle S_n\rangle$ is not finitely presented.}
\end{quote}

\begin{prop} \label{FxF}
If $F$ is a non-abelian free group,
$F\times F$ has property $\Psi$.
\end{prop}

This result is due to C.F.~Miller III
(\cite{cfm-thesis}, \cite{LS} page 194
and \cite{cfm-survey}). 

\begin{proof}
It is
enough to consider $X=\{(a,1), (1,a)\mid a\in A\}$
where $A$ is a basis for $F$. Associated to each finite
group presentation $\P=\ARp$ one has the fibre product
$P=\text{\rm{gp}}\{(a,a),(r,1)\mid a\in A, r\in R\}
\subset F\times F$.
Let $S(\P)$ be (the obvious words representing) the given 
generators of $P$. If $\P$ is trivial then $P=F\times F$.
If $\ARp$ is an infinite group then $P$ is not
finitely presented (see \cite{grun}).

It is well known that there exist recursive sequences
of finite presentations $\P_n = \langle A\mid R_n\rangle$, with
$R_n$ of fixed cardinality, such that each group 
presented is either trivial or infinite 
and there is no algorithm to recognise which
are trivial (see \cite{cfm-thesis}). Define
$S_n= S(\P_n)$.
\end{proof}
 
By applying a modification
of the Rips construction \cite{rips} to
the standard presentation of $F\times F$ one can translate
the undecidability expressed in Proposition \ref{FxF} into the
realm of hyperbolic groups (cf.~\cite{bms} and \cite{BW}).

\begin{prop} \label{Rips}
There exists a  compact negatively curved
2-complex   whose fundamental group  has  
property $\Psi$.
\end{prop}

Recall that a subgroup $H$ of a finitely generated group
$\G$ is said to have an {\em unsolvable membership problem}
if there does not exist an algorithm that takes as input
words $w$ in the generators of $\G$ and decides whether
or not $w\in H$.

\smallskip

\noindent{\em Notation.}
We write $Z_G(g)$ to denote the centralizer in a group $G$
of an element $g\in G$ and $Z_G(S)$ to denote the centralizer of
a subset $S\subset G$.

\begin{lemma} \label{membFF}
If $F$ is a non-abelian free group 
then there exist finitely generated
subgroups $H\subset F\times F$ and elements $h\in H$ such
that $H$ has 
an unsolvable membership problem and  
$Z_{F\times F}(h)\subset H$.
\end{lemma}

\begin{proof} Let $A$ be a basis for $F$ and let $\ARp$
be a finitely presented group whose word problem is
unsolvable. Consider the associated fibre product $P\subset F\times
F$, as in the proof of (\ref{FxF}). Now, given a word
$w$ in the letters $A^{\pm 1}$, ask if $(w,1)\in P$. The
answer is YES  if and only if $w=1$ in  $\ARp$, and the
validity of this equality
cannot be determined algorithmically. Define $H=P$ and
$h=(a,a)$, where $a\in A$.
\end{proof}

Once more, by applying the Rips construction one can translate
the above phenomenon into the realm of hyperbolic groups ---
see \cite{BH} p.488.

\begin{lemma} \label{membRips}
There exists a compact negatively curved
2-complex  whose fundamental group $\G$ has the  following
property: there is a finitely generated
subgroup  $H\subset \G$ and an element $h\in H$ such
that $H$ has 
an unsolvable membership problem and  $Z_{\G}(h)\subset H$.
\end{lemma}

\section{Further Preliminaries}

This section contains four lemmas that we shall need in 
the proofs of
 Theorems \ref{thmA} and \ref{thmB}. The first two
lemmas concern conjugacy and 
centralizers in
amalgamated free products. These are  special cases
of standard results that are proved in a combinatorial
manner in  \cite{LS}, for example; they also admit simple
geometric proofs in the language of Bass-Serre theory.

Consider an amalgamated free product $\G=A\ast_C B$.

\begin{lemma} \label{conj-ab}
For all $a,a'\in A\ssm C$ and all $b,b'\in B\ssm C$,
if $ab$ is conjugate to $a'b'$ in $\G$ then
there exists $c\in C$ such that $c(ab)c^{-1}=a'b'$.
\end{lemma}

\begin{lemma} \label{centACB} If $S\subset A\ssm C$ then
$Z_\G(S)=Z_A(S)$. 
\end{lemma} 

Our third lemma also involves Bass-Serre theory. Recall that
a group is said to have {\em property FA} if every action of the
group on a tree has a fixed point (equivalently, the group
cannot be written as a non-trivial amalgamated free product
or HNN extension).

\begin{prop} \label{ABQ}
Let $A$ be a group 
that has property FA,
let $Q$ be a  group such that there are no
non-trivial homomorphisms $A\to Q$, and let $B$ be
a group that does not contain an isomorphic
copy of $A$. Let $G=(A\ast B)\times Q$, let
$\phi: Q\to B$ be a homomorphism, and let
$\hat Q = \{(\phi(q),q)\mid q\in Q\}\subset G$.
Let $D$ be the amalgamation of
two copies of $G$ along $\hat Q$. 

Then $K_\phi:=\ker\phi\subset\hat Q$
is a characteristic subgroup of $D$. More precisely,
$K_\phi$ is the intersection of
the centralizers of those subgroups $H\subset D$ that
are isomorphic to $A$ and are not strictly contained in
any subgroup isomorphic to $A$. 
\end{prop}

\begin{proof}  
First 
we claim that $H$ must lie in a conjugate of
one of the obvious copies of  $G$ in $D$. 
To see that this is the case, 
consider the action of $D=G\ast_{\hat Q}\-G$
on the Bass-Serre tree associated to the given decomposition:
$H\cong A$ has property
FA and hence is contained in a vertex stabilizer. 

Now, since there are no non-trivial homomorphisms from $A$
to $Q$, if $H\subset G=(A\ast B)\times Q$ is
isomorphic to $A$ then  $H\subset A\ast B$.
Since $A$ has property FA, it is {\em a fortiori}
freely indecomposable, so the Kurosh subgroup
theorem tells us that $H\cong A$ is contained in a conjugate
of one of the vertex groups
$A, B\subset A\ast B$. And since $B$ does not contain a
copy of $A$, only the former possibility can occur.
It follows that if $H$ is maximal among subgroups
isomorphic to $A$  then it must equal a conjugate
of $A$.  Thus  the maximal $H\subset D$ isomorphic
to $A$ are exactly the conjugates of the obvious subgroups
$A\subset G\subset D$ and $\-A\subset \-G\subset D$.

It remains to calculate 
the intersection of the centralizers of the conjugates
of $A$ and $\-A$. It is easy to check that this intersection
contains $\ker\phi\subset\hat Q$. For the reverse inclusion,
note that since $A\cap\hat Q$ is trivial,
the centralizer of $A$ in $D=G\ast_{\hat Q}\-G$ is simply
its centralizer in $G$ (Lemma \ref{centACB}),
 which is $Z(A)\times Q$. Likewise,
if $b\in B$ is non-trivial, then the centralizer of
$A^b=bAb^{-1}$ in $D$ is $Z_G(bAb^{-1})= bZ_G(A)b^{-1}\times Q$.
Thus the centralizer in $D$ of $A\cup A^b$ is
$Q$. Similarly, the centralizer in $D$ of $\-A\cup \-A^{\-b}$ is
$\-Q$. And $Q\cap \-Q = Q\cap\hat Q = \ker\phi$.
\end{proof}

In our proof of Theorem \ref{thmB} we shall need the following
variation on Rapaport's Theorem \cite{rap}.

\begin{lemma} \label{rapa}
Let $Q$ be a group generated by $\{g_1,\dots,g_m\}$ and let 
$F$ be a free group with basis
$\{\s_1,\dots,\s_m,\tau_1,\dots,\tau_m\}$.
 Let $\phi: F\to Q$ be an 
epimorphism such that $\phi(\tau_i)=1$ for $i=1,\dots,m$. Then
there exists an automorphism $\Phi: F\to F$ such that
$\phi\circ\Phi (\s_i) =1$ and $\phi\circ\Phi(\tau_i)=g_i$ for
$i=1,\dots,m$.
\end{lemma}

\begin{proof} Since $\phi$ is onto, for each $i=1,\dots,m$
there exists a word $u_i$ in the letters $\s_j$ and
their inverses such that $g_i=\phi(u_i)$ in $Q$. And
since the $g_j$ generate $Q$, for
each $i=1,\dots,m$ there exists a word $v_i$ in the
letters $g_j$ and their inverses such that $\phi(\s_i) = v_i$ in $Q$.

Let $\Phi_1: F_{2m}\to F_{2m}$ be the automorphism that sends
each $\s_i$ to itself and sends $\tau_i$ to $\tau_i':=\tau_iu_i$
for $i=1,\dots,m$.

Let $V_i$ be the word obtained from $v_i$ by substituting the
word $\tau_j'$ for each occurrence of the symbol $g_j$, for 
$j=1,\dots,m$.

Let $\Phi_2: F_{2m}\to F_{2m}$ be the automorphism that sends
each $\tau_i'$ to itself and sends $\s_i$ to $\s_i':=\s_iV_i^{-1}$
for $i=1,\dots,m$.

(To see that  $\Phi_1$ and $\Phi_2$ are automorphisms,
note that each is the composition of right Nielsen moves.)

Define $\Phi = \Phi_2\circ\Phi_1$.
\end{proof} 
  
\section{Unsolvable Conjugacy Problems}

In this section we shall prove Theorem \ref{thmA}.

\begin{definition} Given a group $G$ and a subgroup
$H\subset G$, we write $\Delta(G;H)$ to denote
the {\em double of $G$ along $H$}, i.e. the free product
of two copies of $G$, written $G$ and $\-G$, with the
subgroups $H$ and $\-H$ amalgamated by the identification
$h=\-h\ \forall h\in H$.
\end{definition}

\begin{thm} \label{basic}
Let $N$ and $F$ be finitely generated groups and 
let $\phi: F\to N$ be a homomorphism. Let
$\hat F = \{(\phi(x),x)\mid x\in F\}\subset N\times F$.

If the membership problem for 
$\im\phi$ is unsolvable and there exists $a_0\in N$ such that
$Z_N(a_0)
\subset \im\phi$, then 
$D=\Delta(N\times F; \hat F)$ has an unsolvable conjugacy problem.
\end{thm}

\begin{proof} Let $I=\im\phi$. We claim that our hypotheses
imply that the following problem is algorithmically unsolvable:

\smallskip

\begin{framebox}{
Given\footnote{Here, as usual, 
``given $b\in N$" should be interpreted as meaning
that one is given a word in a fixed 
finite generating set for $N$.} $b\in N$, 
determine whether $\exists h\in I$ with 
$ha_0h^{-1}=ba_0b^{-1}$.}
\end{framebox}

\smallskip

To see this, first note
that $ha_0h^{-1}=ba_0b^{-1}$ if and only $h^{-1}b\in
Z_D(a_0)$. Since $h$ and $Z_D(a_0)$ are assumed to
lie in $I$, this means that such an $h$ exists
if and only if $b\in I$. And we are assuming that there is
no algorithm to decide whether $b\in I$.

We adopt the notation $\hat x = (\phi(x),x)$ for elements of
$\hat F$. 
Identifying $N$ and $F$ with
their images in $N\times F$ we have
$$
\hat x n \hat x^{-1} = \phi(x)\, n\, \phi(x)^{-1}\ \forall n\in N\
\forall x\in F.
$$

To see that the conjugacy problem in $D$ is unsolvable  we
consider the following sub-problem:

\smallskip

\begin{framebox}{
given $b\in N$, decide if $(ba_0b^{-1})(\-{ba_0b^{-1}})$
is conjugate to $a_0\-a_0$ in $D$.}
\end{framebox}
\smallskip

Lemma \ref{conj-ab} tells us that 
$(ba_0b^{-1})(\-{ba_0b^{-1}})$
is conjugate to $a_0\-a_0$ in $D=(N\times F)\ast_{\hat F}
(N\times F)$ if and only if there exists $\hat x\in\hat F$
conjugating one element to the other, i.e.
$$
\hat x a_0 \-a_0 \hat x^{-1} =
(ba_0b^{-1})(\-{ba_0b^{-1}}).
$$
And noting that $\hat x = \-{\hat x}$,  we have

$$
\hat x a_0 \-a_0 \hat x^{-1} = (\hat x a_0 
\hat x^{-1})\-{(\hat x a_0 \hat x^{-1})} =
(\phi(x) a_0 \phi(x)^{-1})\-{(\phi(x) a_0 \phi(x)^{-1})}.
$$

Thus $(ba_0b^{-1})(\-{ba_0b^{-1}})$
is conjugate to $a_0\-a_0$ in $D$ if and only if
there exists an element $h\in I$ such that
$(ba_0b^{-1})(\-{ba_0b^{-1}})= (ha_0h^{-1})(\-{ha_0h^{-1}}).$
Moreover, since $N, \-N\subset D$ do not intersect the amalgamated subgroup
$\hat F$, the subgroup they generate is the free product
$N\ast\-N$. Therefore
$(ba_0b^{-1})(\-{ba_0b^{-1}})= (ha_0h^{-1})(\-{ha_0h^{-1}})$ if and
only if $ba_0b^{-1} = ha_0h^{-1}$.

Thus we have shown that, given $b\in N$, 
the elements $(ba_0b^{-1})(\-{ba_0b^{-1}})$ and $a_0\-a_0$
are conjugate in $D$ if and only if there exists
$h\in I$ such that $ba_0b^{-1} = ha_0h^{-1}$. In the
first paragraph of this proof we explained that there
is no algorithm to decide whether such an $h$ exists.
\end{proof}

We use the term {\em semihyperbolic} in the sense of \cite{AB}.
In particular, biautomatic groups are semihyperbolic, as are
the fundamental groups of compact non-positively curved 
(orbi)spaces.

As a special
case of the main construction in \cite{mb-cmh} we have the
following:

\begin{thm} \label{cmh}
Let $N$ and $F$ be groups,
let $\phi:F\to N$ be a homomorphism, and
let $\hat F = \{(\phi(x),x)\mid x\in F\}$.

If $N$ is semihyperbolic and $F$ is combable, then 
$\Delta(N\times F; \hat F)$ is combable.
\end{thm}

\begin{para}{\bf The Proof of Theorem \ref{thmA}}\end{para}

Lemmas \ref{membFF} and \ref{membRips} provide examples of 
semihyperbolic groups $N$ and finitely generated subgroups
$H\subseteq N$ such that the membership problem for $H$ is
unsolvable, and there exists $h\in H$ such that $Z_N(h)\subset H$.
Fix a surjection $\phi:F\to H$, where $F$ is 
a finitely
generated free group. Theorem \ref{cmh} tells us that
$D=\Delta(N\times F; \hat F)$ is combable and Theorem \ref{basic}
tells us that $D$ has an unsolvable conjugacy problem. 
\hfill$\square$

\begin{remark} The 
results in \cite{mb-cmh} show that each of the examples $D$
constructed in the above proof has a cubic Dehn function and is 
{\rm{Ind}}-combable, where {\rm{Ind}} is the class of indexed
languages. At the time of writing, I do not know of
similar examples that are {\rm{Reg}}-combable, i.e.
automatic.  
\end{remark}

\section{Unsolvable Isomorphism Problems}

Our purpose in this section is to prove Theorem \ref{thmB}. We
shall achieve this by applying the 
following general criterion to the groups from
Section 1 that satisfy Property $\Psi$.  

\begin{thm} \label{thmB'}
Let $A$ be a group with property FA and
let $B$ be a finitely presented group that does not contain a 
subgroup isomorphic to $A$. 
If $B$ has property $\Psi$ and $F$ is a finitely generated
subgroup of sufficiently large rank $2m$,
then there is a recursive
sequence of subsets $S_n=\{s_{n,1},\dots,s_{n,m}\}
\subset B$
such there is no algorithm to determine which of the groups
$$
\G_n= \Delta((A\ast B)\times F; \S_n)
$$
are isomorphic, where $\S_n \subset B\times F$
is the subgroup $\langle(s_{n,i},\sigma_i),\, (1,\tau_i)
; i=1,\dots,m\rangle$, with
$\{\sigma_1,\dots,\sigma_m,\tau_1,\dots,\tau_m\}$ a 
fixed basis for 
$F$.
\end{thm}
\smallskip

Before turning to the proof of this result we note that
the groups $\G_n$  are obtained in a natural way 
from a  fixed free group by adding a fixed number of relations.
Specifically, if we fix  a presentation  $\langle
\Cal X\mid \Cal R\rangle$ for $(A\ast B)\times F$, with
$\{\s_1,\dots,\s_m,\tau_1,\dots,\tau_m\}\subset\Cal X$
a basis for $F$, 
then 
$$
\G_n = \langle \Cal X,
\-{\Cal X} \mid \Cal R, \-{\Cal R};
\tau_i=\-\tau_i,\, \s_is_{n,i}=\-{\s_is_{n,i}}
\text{ for } i=1,\dots,m\rangle.
$$

\begin{proof}  
$F$ is the free group with basis $\{\s_1,\dots,\s_m,
\tau_1,\dots,\tau_m\}$.
The subgroup $\S_n\subset B\times F$ is  
$\{(\phi_n(x),x) \mid x\in F\}$, where 
$\phi_n : F\to B$ is the homomorphism that sends each
$\s_i$ to $s_{n,i}\in S_n$ and sends each $\tau_i$ to the
identity. The hypothesis that $B$ has property $\Psi$ means
that we may assume that there is no algorithm to determine
which of the subgroups $\langle S_n\rangle$ are equal
to $B$; we may also assume that those which
are not equal are not finitely presentable. We
proceed with these  assumptions in place.
\smallskip
\noindent{\em Scheme of Proof:} In Step 1 we prove that 
 if
$\langle S_n\rangle\subset B$ is finitely presentable
and $\langle S_{n'}\rangle\subset B$ is not, then
$\G_n\not\cong\G_{n'}$. In Step 2 we argue that if
$\langle S_n\rangle\subset B$ and 
$\langle S_{n'}\rangle\subset B$  are both 
finitely presentable (hence equal to $B$) then 
$\G_n\cong\G_{n'}$.
Since there is
no algorithm to determine which of the subgroups 
$\langle S_n\rangle\subset B$ are finitely presentable,
this will complete the proof.
\smallskip

\noindent{\em Step 1:}
Proposition \ref{ABQ} tells us that any
isomorphism from $\G_n$  
to $\G_{n'}$ must map
$\ker\phi_n\subset\S_n$ to $\ker\phi_{n'}$ isomorphically. 
(We may apply Proposition \ref{ABQ} because no group
with property FA admits a non-trivial homomorphism to
a free group.) Thus we shall be done if we can prove that
$\ker\phi_n$ is the normal
closure of a finite subset of $\G_n$ if and only if
$\langle S_n\rangle$ is finitely presentable.

$\phi_n$ is an epimorphism from the finitely generated free
group $F$ to the group $\langle S_n\rangle$. 
Thus $\ker\phi_n$
is the normal
closure of a finite subset of $F$ if and only if 
$\langle S_n\rangle$ is finitely presentable. (The
 finite presentability of a group is independent of
the choice of finite generating set.)

There is a natural retraction of $\G_n$ onto $F$: first,
writing $G=(A\ast B)\times F$,
retract  $\G_n=G\ast_{\Sigma_n}\-G$ onto 
$G$ by sending both  $g\in G$ and  $\-g\in\-G$ to
$g$ for every $g\in G$; now project onto the
second direct summand in $G$. 

If a subgroup
$H$ is a retract of a group $\G$, then  the normal closure
in $H$ of any subset of $U\subset H$ is the same as the
image in $H$ (under retraction) of the normal closure
of $U$ in $\G$. Thus, combining the 
conclusions of the previous two paragraphs,
we see that $\ker\phi_n$
is the normal
closure of a finite subset of $\G_n$ if and only if 
$\langle S_n\rangle$ is finitely presentable.

\smallskip
\noindent{\em Step 2:} We
now assume that  $\langle S_n\rangle\subset B$ and 
$\langle S_{n'}\rangle\subset B$  are both 
finitely presentable. In this case $\phi_n$ and
$\phi_{n'}$ have the same image, namely $B$. Thus
we may apply Lemma \ref{rapa}. 
This lemma provides us
with an automorphism $\Phi: F\to F$ such that
$\phi_n\circ\Phi (\s_i) =1$ and $\phi_n\circ\Phi(\tau_i)=s_{n',i}$ for
$i=1,\dots,m$.

Consider
the automorphism $\alpha$ of $G=(A\ast B)\times F$ that
is the identity on $A\ast B$ and is defined on $F$ by
$\alpha(\s_i)=\Phi(\tau_i)$ and $\alpha(\tau_i)=\Phi(\sigma_i)$.
We claim  that $\alpha$
sends $\S_{n'}$ isomorphically to $\S_{n}$. Since 
these are free groups of equal rank, it suffices
to check that  $\alpha(\S_{n'})\subset \S_{n}$.

By definition, for 
$i=1,\dots,n$ we have $\alpha(s_{n',i},\sigma_i)
= (s_{n',i}, \Phi(\tau_i))$, which is in 
$\S_n=\{(\phi(x),x)\mid x\in F\}$  because
$\phi_n\circ\Phi(\tau_i) =  s_{n',i}$.  

Because $\alpha$ maps $\S_{n'}$ isomorphically onto
$\S_{n}$, we obtain an isomorphism 
from $\G_{n'}=G\ast_{\S_{n'}}\-G$ to
$\G_{n}$ by defining $g\mapsto 
\alpha(g)$ and $\-g\mapsto\-{\alpha(g)}$.

This completes the
proof of Theorem \ref{thmB'}.
\end{proof}

\subsection{The Proof of Theorem \ref{thmB}}
 
There are many semihyperbolic groups that have property FA, for
example finite groups, mapping class groups, higher-rank lattices,
and hyperbolic groups that have property (T).
To be definite, we take $A$ to be a finite group and
we take $B$ to be a 
torsion-free semihyperbolic group that has property $\Psi$ (examples
are discussed below). By applying the construction of 
Theorem \ref{thmB'} to $G:=(A\ast B)\times F$ we obtain 
a recursive sequence of groups $\G_n = \Delta(G;{\S_n})$ for
which the isomorphism problem is unsolvable. According
to Theorem \ref{cmh}, the groups  $\G_n $ are all
combable. 
\hfill$\square$

\begin{example} Let $F$ be a finitely generated free group
of rank at least 2. The arguments in Section 1 yield 
explicit sequences of subsets $S_n\subset F\times F$ 
of cardinality $m$ enjoying  Property $\Psi$.
 Let $F_m$ be a free group of rank $2m$.
Let $G=[\Bbb Z_2\ast (F\times F)]
\times F_m$. The construction of Theorem \ref{thmB'}
shows that there is a recursive sequence
of  free subgroups $\S_n\subset G$ (each a 
retract of $G$ in fact) such that there is no
algorithm to decide isomorphism among the combable
groups $\G_n=\Delta(G;\S_n)$.
\end{example}

\begin{example} The arguments in Section 1 also yield 
torsion-free
hyperbolic groups $\G$ with
explicit sequences of subsets $S_n\subset\G$ 
of cardinality $m$ enjoying  Property $\Psi$.
As above, in $G=(\Bbb Z_2\ast\G)\times F$ one
finds a recursive sequence
of  free subgroups $\S_n\subset G$ (each a 
retract of $G$) such that there is no
algorithm to decide isomorphism among the combable
groups $\G_n=\Delta(G;\S_n)$. 
\end{example}

\begin{remark} If one wishes to obtain torsion-free examples, one
can replace $\Bbb Z_2$ in the above constructions
by a torsion-free hyperbolic
group with Property (T).
\end{remark}

\begin{remark} The results in \cite{mb-cmh} show that all of
the groups in the   sequences $(\G_n)$ considered above
satisfy a cubic Dehn isoperimetric inequality and are
{\rm{Ind}}-combable, where {\rm{Ind}} is the class of indexed
languages. At the time of writing, I do not know of such
sequences in which all of the groups are {\rm{Reg}}-combable, i.e.
automatic.
\end{remark}

\enddocument